\documentclass [12pt,a4paper,reqno]{amsart}
\input amssymb.sty
\usepackage{amsmath,amsthm}
\usepackage{amsfonts}
\allowdisplaybreaks

\newcommand{\be}{\begin{equation}}
\newcommand{\ef}{\end{equation}}
\chardef\bslash=`\\ 

\newtheorem{thm}{Theorem}[section]
\newtheorem*{thm*}{Theorem}

\newtheorem{lem}[thm]{Lemma}
\newtheorem{prop}[thm]{Proposition}

\theoremstyle{definition}

\newtheorem*{remark*}{Remarks}
\newtheorem*{defn*}{Definition}

\theoremstyle{remark}

\numberwithin{equation}{section}

\newcommand{\G}{\Gamma}
\newcommand{\wt}{\widetilde}
\newcommand{\wh}{\widehat}
\newcommand{\fc}{\frac}
\newcommand{\bk}{\bigskip}
\newcommand{\iy}{\infty}
 \renewcommand{\sectionmark}[1]{}

\renewcommand{\Re}{\operatorname{Re}}

\newcommand{\hol} {holomorphic}
\newcommand{\qc} {quasiconformal}
\newcommand{\sh} {subharmonic}
\newcommand{\psh} {plurisubharmonic}

\newcommand{\ve}{\varepsilon}
\newcommand{\e}{\epsilon}

\newcommand{\Ko} {Kobayashi}
\newcommand{\Ca} {Carath\'{e}odory}

 \usepackage{amsfonts}
\newcommand{\field}[1]{\mathbb{#1}}
\newcommand{\g}{\gamma}

\newcommand{\dl}{\delta}
\newcommand{\D}{\Delta}
\newcommand{\om}{\omega}
\newcommand{\z}{\zeta}
\newcommand{\ov}{\overline}
\newcommand{\vp}{\varphi}

\newcommand{\C}{\field{C}}

\newcommand{\B}{\mathbf{B}}
\newcommand{\T}{\mathbf{T}}

\newcommand{\Hol}{\operatorname{Hol}}

\newcommand{\Om} {\Omega}

\newcommand{\x} {\mathbf x}

\renewcommand{\a} {\alpha}

\newcommand{\ld}{\lambda}
\newcommand{\kp}{\kappa}

\newcommand{\Pot}{\operatorname{Pot}}
\newcommand{\hyp}{\operatorname{hyp}}
\newcommand{\Mob}{\operatorname{Mob}}

\begin{document}

\title{Hyperbolic distances, nonvanishing holomorphic functions
and Krzyz's conjecture}
\author{Samuel L. Krushkal}

\begin{abstract} The goal of this paper is to prove the
conjecture of Krzyz posed in 1968 that for nonvanishing
\hol \ functions  $f(z) = c_0 + c_1 z + \dots$ in the
unit disk with $|f(z)| \le 1$, we have the sharp bound
$|c_n| \le 2/e$ for all $n \ge 1$,
with equality only for the function
$f(z) = \exp [(z^n - 1)/(z^n + 1)]$
and its rotations.
The problem was considered by many researchers, but only partial results have been established.
The desired estimate has been proved only for $n \le 5$.

Our approach is completely different and relies on
complex geometry and pluripotential features of convex
domains in complex Banach spaces.
\end{abstract}

\date{\today\hskip4mm({krzyz.tex})}

\maketitle

\bigskip

{\small {\textbf {2000 Mathematics Subject Classification:}
Primary: 30C50, 30C55, 32Q45; Secondary 30F45, 31C10, 32U35}

\medskip

\textbf{Key words and phrases:} Nonvaninshing \hol \ function,
Krzyz's conjecture, \Ko \ metric, \Ca \ metric, complex geodesic,
\psh \ function, pluricomplex Green function, convex domain}

\bigskip

\markboth{Samuel L. Krushkal}{Hyperbolic distances, Krzyz's conjecture}
\pagestyle{headings}

\bk

\section{Krzyz's conjecture. Main theorem}.

Nonvanishing \hol \ functions $f(z) = c_0 + c_1 z + ...$ on the unit
disk $\D = \{z : |z| < 1\}$ (i.e., such that $f(z) \ne 0$ in $\D$)
 form the normal families admitting
certain invariance properties, for example, the invariance
under action of the
M\"{o}bius group of conformal self-maps of $\D$, complex homogeneity, etc.
One of the most interesting examples of such families is the set
$\mathcal B_1 \subset H^\iy$ of \hol \ maps of $\D$ into
the punctured disk $\D_{*} = \D \setminus \{0\}$.

Compactness of $\mathcal B_1$ in topology of locally uniform convergence
on $\D$ implies the existence for each $n \ge 1$ the
extremal functions $f_0$ maximizing $|c_n(f)|$ on $\mathcal B_1$.
Such functions are nonconstant and must satisfy
$|f(e^{i \theta})| = 1$ for almost all $\theta \in [0, 2 \pi]$.

The problem of estimating coefficients on $\mathcal B_1$ was posed by
Krzyz \cite{Kz} in 1968. He conjectured that for all $n \ge 1$,
\be\label{1.1}
|c_n| \le \fc{2}{e},
\end{equation}
with equality only for the function
\be\label{1.2}
\kappa_0(z) := \exp \Bigl(\frac{z -1}{z + 1} \Bigr) =
\frac{1}{e} + \frac{2}{e} z - \frac{2}{3e} z^3 + ... \ .
\end{equation}
and its rotations
$\epsilon_1 \kappa_0(\epsilon_2 z)$ with $|\epsilon_1| = |\epsilon_2| = 1$.
Note that (1.2) provides a \hol \ universal covering map
$\D \to \D_{*}$ with $f(0) = 1/e$.

This fascinating and extremely interesting problem has been investigated
by a large number of mathematicians, however it still remains open.
The estimate (1.1) has been proved
only for $n \le 5$ (see \cite{HSZ}, \cite{PS}, \cite{Sa}, \cite{Sz},
\cite{Ta}).

The best uniform estimate for all $n$ given by Horowitz \cite{Ho} is
$$
|c_n| \le 1 - \fc{1}{3 \pi} + \fc{4}{\pi}
\sin \Bigl(\frac{1}{12} \Bigr) = 0.999...
$$
(while $2/e = 0.7357...)$; it was somewhat improved later.
For a more complete history of this problem we refer e.g., to \cite{Ba},
\cite{HSZ}, \cite{LS}, \cite{Sz}.

Our goal is to prove that Krzyz's conjecture is true for all $n \ge 1$:

\begin{thm} For every $fz) = c_0 + c_1 z + ... \in \mathcal B_1$
and $n \ge 1$, we have the sharp bound (1.1), and the equality occurs only
for the function (1.2) and its rotations.
\end{thm}

Our approach is completely different and relies on complex geometry
and pluripotential features of convex domains in complex Banach spaces.
The underlying idea of the proof is in fact the same as for Zalcman's
conjecture applied in \cite{Kr4} (and earlier in \cite{Kr3}).
It uses also the important fact that the function (1.2)
generates the complex geodesics in a domain formed by nonvanishing
functions on the closed unit disk.
Certain results obtained in the proof of the main theorem have
independent interest.
Let us mention also that the proof essentially involves certain
specific features of $H^\iy$.

\section{Preliminaries: Hyperbolic metrics on convex Banach domains}

We first present briefly the basic results on
properties of the \Ko \ and \Ca \ metrics and on complex geodesics
on convex domains in complex Banach spaces, which underly the
proof of Theorem 1.1.

\subsection{Equality of metrics}
Let $D$ be a complex Banach manifold modelled by a Banach space $X$.
The  {\bf \Ko \ metric} $d_D$ on $D$ is the
largest pseudometric $d$ on $D$ that does not get increased by \hol \ maps
$h: \ \D \to D$ so that for any two points $\x_1, \ \x_2 \in D$,
we have
$$
d_D(\x_1, \x_2) \le \inf \{d_\D(0,t): \ h(0) = \x_1, \ h(t) = \x_2\},
$$
where $d_\D$ is the {\bf hyperbolic Poincar\'{e} metric} on $\D$ of Gaussian
curvature $- 4$, with the differential form
\be\label{2.1}
ds = \ld_{\hyp}(z)|dz| : = |dz|/(1 - |z|^2).
\end{equation}
The {\bf \Ca} \ distance between $\x_1$ and $\x_2$ in $D$ is
$$
c_D(\x_1, \x_2) = \sup d_\D(f(\x_1), f(\x_2)),
$$
where the supremum is taken over all \hol \ maps $f: \ \D \to X$.

The corresponding {\bf differential} (infinitesimal) forms of the \Ko \
and \Ca \ metrics are defined for the points
$(\x, v)$ in the tangent bundle $T D$ of $D$, respectively, by
$$
\begin{aligned}
\mathcal K_D(\x, v) &= \inf \{r: \ r > 0, \ \exists \ h \in \Hol(\D, D), \
h(0) = \x, \ d h(0) r = v\},  \\
\mathcal C_D(\x, v) &= \sup \{|d f(\x) v|: \ f \in \Hol(\T, \D), \ f(\x) = 0\},
\end{aligned}
$$
where $\Hol(X, Y)$ denotes the collection of \hol \ maps of a complex
manifold $X$ into $Y$. Note that in the general case,
\be\label{2.2}
\limsup\limits_{t\to 0, \ t \ne 0}
\fc{d_D(\x, \x + t \mathbf v)}{|t|} \le \mathcal K_D(\x, \mathbf v).
\end{equation}
For general properties of invariant metrics we refer to \cite{Di},
\cite{Ko}.
A remarkable fact is:

\begin{prop} If $D$ is a convex domain in complex Banach space, then
\be\label{2.3}
d_D(\x_1, \x_2) = c_D (\x_1, \x_2) =
\inf \{d_\D (h^{-1}(\x_1), h^{-1}(\x_2)): \ h \in \Hol(\D, D)\}
\end{equation}
and
\be\label{2.4}
\mathcal K_D(\x, v) = \mathcal C_D(\x, v) \quad \text{for all} \ \
(\x, v) \in T(D).
\end{equation}
In particular, both infinitesimal and global pseudo-distances are
logarithmically \psh \ on $D$.
\end{prop}

In the case of a bounded domain $D$, both $d_D$ and $c_D$ are distances
(metrics), which means that these geometric quantities separate the
points in $D$.

The equality of global pseudo-distances on convex domains
in $\C^n$ and their representations by (2.3) were established by
Lempert \cite{Le}; the coincidence of the infinitesimal
metrics for such domains was proved by Royden and Wong \cite{RW}.
These results were extended to convex domains in infinite dimensional
Banach spaces in Dineen-Timoney-Vigu\'{e} \cite{DTV}.

\subsection{Pluripotential and curvature properties}
Proposition 2.1 is rich in corollaries. We shall use several of them.

First recall that the {\bf pluricomplex Green function} $g_D(x, y)$
of a domain $D \subset X$ with pole $\mathbf y$ is defined by
\be\label{2.5}
g_D(\x, \mathbf y) = \sup u_{\mathbf y}(\x)
\quad (\x, \mathbf y \in D)
\end{equation}
and following upper regularization
$$
v^*(\x) = \limsup_{x^\prime \to \x} v(\x^\prime).
$$
The supremum  in (2.5) is taken over all \psh \ functions
$u_{\mathbf y}(\x): \ D \to [-\iy, 0)$ such that
$$
u_{\mathbf y}(\x) = \log \|\x - \mathbf y\|_X + O(1)
$$
in a neighborhood of the pole $\mathbf y$; here $\|\cdot\|_X$ denotes
the norm
on $X$, and the remainder term $O(1)$ is bounded from above (cf. e.g.,
[Di]).
The Green function $g_D(\x, \mathbf y)$ is a maximal \psh \ function on
$D \setminus \{\mathbf y\}$ (unless it is not identically $- \iy$).
Proposition 2.1 implies

\begin{prop} If $D$ is a convex domain in a complex Banach space, then
\be\label{2.6}
g_D(\x, \mathbf y) = \log \tanh d_D(\x, \mathbf y)
= \log \tanh c_D(\x, \mathbf y)
\end{equation}
for all $\x, \mathbf y \in D$.
\end{prop}

The next corollary concerns the curvature properties.
There are several generalizations of the smooth
Gaussian curvature.

The generalized Gaussian curvature $\kappa_\ld$ of an upper
semicontinuous Finsler (semi)metric $ds = \ld(t) |dt|$ in a domain
$\Om \subset \C$ is defined by
\be\label{2.7}
\kappa_\ld (t) = - \fc{\mathbf{\D} \log \ld(t)}{\ld(t)^2},
\end{equation}
where $\mathbf{\D}$ is the {\bf generalized Laplacian} defined by
\be\label{2.8}
\mathbf{\D} \ld(t) = 4 \liminf\limits_{r \to 0} \frac{1}{r^2}
\Big\{ \frac{1}{2 \pi} \int_0^{2\pi} \ld(t + re^{i \theta}) d
\theta - \ld(t) \Big\}
\end{equation}
(provided that $- \iy \le \ld(t) < \iy$).
It is well-known that an upper semicontinuous function $u$ is \sh \
on its domain $D \subset \C$ if and only if
$\mathbf \D u(t) \ge 0$ on its domain $D \subset \C$;
hence, at the points $t_0$ of local
maximuma of $\ld$ with $\ld(t_0) > - \iy$, we have $\mathbf{\D}
\ld(t_0) \le 0$.
Note that for $C^2$ functions, $\mathbf \D$ coincides with
the usual Laplacian
$4 \partial^2/\partial z \partial \ov z$, and its non-negativity
immediately follows
from the mean value inequality; for arbitrary \sh \ functions,
this is obtained by a standard approximation.

The {\bf sectional \hol \ curvature} of a Finsler metric on a
complex Banach manifold $D$ is defined in a similar way as the
supremum of the curvatures (2.7) over appropriate collections of
\hol \ maps from the disk into $D$ for a given tangent direction
in the image.

The \hol \ curvature of the \Ko \ metric $\mathcal K_D(\x, v)$ of any
complete hyperbolic manifold $D$ satisfies $\kp_{\mathcal K_D}(\x, v) \ge
- 4$ at all points $(\x, v)$ of the tangent bundle $T D$
of $D$, and for the \Ca \ metric $\mathcal C_D$ we have
$\kp_{\mathcal C_D}(\x, v) \le - 4$ (see e.g., [Di]). Consequently, at
each point, where these metrics are equal, we have the equality
\be\label{2.9}
\kp_{\mathcal K_D}(\x, v) = \kp_{\mathcal C_D}(\x, v) = - 4.
\end{equation}
By Proposition 2.1, this holds for all convex domains $D$.

\bk
It follows from (2.7) that a conformal Finsler
metric $ds = \ld(z) |dz|$ with $\ld(z) \ge 0$ of
generalized Gaussian curvature at most $- K, \ K > 0$,
satisfy the inequality
\be\label{2.10}
\D \log \ld \ge K \ld^2,
\end{equation}
where $\D$ is the generalized Laplacian (2.8).
We shall use its integral generalization due to Royden \cite{Ro}.

A conformal metric $\ld(z) |dz|$ in a domain $G$ on $\C$
(more generally, on a Riemann surface) has the curvature less
than or equal to $K$
{\bf in the supporting sense} if for each $K^\prime > K$
and each $z_0$ with
$\ld(z_0) > 0$, there is a $C^2$-smooth supporting metric $\wt \ld$
for $\ld$ at $z_0$ (i.e., such that $\wt \ld(z_0) = \ld(z_0)$
and $\wt \ld(z) \le \ld(z)$ in a neighborhood of $z_0$) with
$\kappa_{\wt \ld}(z_0) \le K^\prime$ (cf. \cite{Ah}, \cite{He}).

A metric $\ld$ has curvature at most $K$ {\bf in the potential sense}
at $z_0$ if there is a disk $U$ about $z_0$ in which the function
$$
\log \ld + K \Pot_U(\ld^2),
$$
where $\Pot_U$ denotes the logarithmic potential
$$
\Pot_U h = \fc{1}{2 \pi} \int\limits_U h(\z) \log |\z - z| d \xi d \eta
\quad (\z = \xi + i \eta),
$$
is \sh. One can replace $U$ by any open subset $V \subset U$, because the
function $\Pot_U(\ld^2) - \Pot_V(\ld^2)$ is harmonic on $U$.
Note that having curvature at most $K$ in the potential sense is
equivalent to $\ld$ satisfying (2.10) in the sense of distributions.

\begin{lem} \cite{Ro} If a conformal metric has curvature at
most $K$ in the supporting
sense, then it has curvature at most $K$ in the potential sense.
\end{lem}

\subsection{Existence of complex geodesics}
Let $D$ be a Banach domain endowed with a pseudo-distance $\rho$.
Following Vesentini (see e.g., \cite{Ve}), a \hol \ map
$h: \ \D \to D$ is called {\bf complex $\rho$-geodesic} if
there exist $t_1 \ne t_2$ in $\D$ such that
$$
d_\D(t_1, t_2) = \rho(h(t_1), h(t_2));
$$
one says also that the points $h(t_1)$ and $h(t_2)$ can be
joined by a complex $\rho$-geodesic.

If $h$ is a complex $c_D$-geodesic, then it is also
$d_D$-geodesic, and vice versa, and then the equality (2.3)
holds for all points of the disk $h(\D)$.

It is important to have the conditions for domains ensuring
the existence and uniqueness of complex geodesics. Certain
conditions, which will be used here, are given in \cite{Di},
\cite{DTV}.

Recall that a Banach space $X$ is called the {\bf dual} of
a Banach space $Y$ if $X = Y^\prime$, that is, $X$ is the
space of bounded linear functionals $x(y) = <x, y>$ on $Y$.
Then $Y$ is called {\bf predual} of $X$.
The {\bf weak$^*$ topology} on $X$
determined by $Y, \ \sigma(X, Y)$, is the topology of pointwise
convergence on points of $Y$, i.e., $x_n \in X \to x \in X$
in $\sigma(X, Y)$ as $n \to \iy$ if and only if
$x_n(y) \to x(y)$ for all $y \in Y$.

If $X$ has a predual $Y$ then the closure $\ov X_1$ of its open
unit ball in $\sigma(X, Y)$ is compact.

\begin{prop} \cite{Di}, \cite{DTV} Let $D$ be a bounded convex
domain
in a complex Banach space $X$ with predual $Y$. If the closure
of $D$ is $\sigma(X, Y)$-compact, then every distinct pair
of points in $D$ can be joined by a complex $c_D$-geodesic.
\end{prop}

This proposition also has its differential counterpart which
provides that under the same assumptions, for any point
$\x \in D$ and any nonzero vector $v \in X$, there exists
at least one complex geodesic $h: \ \D \to D$ such that
$h(0) = \x$ and $h^\prime(0)$ is colinear to $v$ (cf. \cite{DTV}).

Note that along a complex geodesic in $D$, the relation (2.2)
reduces to the equality.

\section{Proof of theorem 1.1}

We prove the main Theorem in several stages;
each stage is of independent interest.

\medskip
\noindent
{\bf 1. Open domain of nonvanishing functions and its
\hol \ embedding}.

\medskip
\noindent
{\em (a)} \ Consider the subsets of $\mathcal B$ defined by
\be\label{3.1}
\mathcal B_r = \{f \in H^\iy(\D_{1/r}): \ f(z) \ne 0 \ \
\text{on the disk} \ \ \D_{1/r} = \{|z| < 1/r\}, \quad
0 < r < 1.
\end{equation}
Note that $\mathcal B_r \subset \mathcal B_{r^\prime}$ if
$r <  r^\prime$.

\begin{lem} Each point of the union
\be\label{3.2}
\mathcal B^0 = \bigcup_r \mathcal B_r
\end{equation}
has a neighborhood in $H^\iy(\D)$, which entirely belongs
to $\mathcal B^0$. Hence, $\mathcal B^0$
is a domain in the space $H^\iy(\D)$.
\end{lem}

\medskip
\noindent
{\bf Proof}. To establish the openness of the union (3.2),
it suffices to show that every function
$f \in \mathcal B_r$
has a neighborhood $U(f, \epsilon(r))$ in $H^\iy(\D)$,
which contains only nonvanishing functions on $\D$.
The connectedness of $\mathcal B^0$ follows from widening
the sets $\mathcal B_r$ when $r$ increases.

Assume the contrary. Then there exist a function
$f_0 \in \mathcal B_r$ and the sequences of functions
$f_n \in H^\iy(\D)$ convergent to $f_0$,
\be\label{3.3}
\lim\limits_{n\to \iy} \|f_n - f\|_{H^\iy(\D)} = 0
\end{equation}
and of points $z_n \in \D$ convergent to $z_0, \ |z_0| \le 1$
such that
$f_n (z_n) = 0 \ (n = 1, 2, \dots)$.

In the case $|z_0| < 1$ we immediately reach a contradiction,
because then the uniform convergence of $f_n$ on
compact sets in $\D$ implies $f_0(z_0) = 0$, which is
impossible.

The case $|z_0| = 1$ requires other arguments.
Since $f_0$ is \hol \ and does not vanish on
the closed disk $\ov \D$,
$$
\min_{|z| \le 1} |f_0(z)| = a > 0.
$$
Hence, for each $z_n$,
$$
|f_n(z_n) - f_0(z_n)| = |f_0(z_n)| \ge a,
$$
and by continuity, there exists a neighborhood
$\D(z_n, \dl_n) = \{|z - z_n < \dl_n\}$ of $z_n$ in $\D$,
in which $|f_n(z) - f_0(z)| \ge a/2$ for all $z$.
This implies
$$
 \|f_n - f_0\|_{H^\iy(\D)} \ge \max_{\D(z_n, \dl_n)} |f_n(z) - f_0(z)|
 \ge \fc{a}{2}.
 $$
This inequality must hold for all $n$, contradicting (3.3).
Lemma follows.

\medskip
\noindent
{\bf Remark}. This lemma does not contradict to existence
of a sequence $\{f_n\} \in H^\iy$, which contains the functions
vanishing in $\D$ and is convergent to $f_0 \in \mathcal B^0$
only uniformly on compact sets in $\D$.

\bigskip
Now put
\be\label{3.4}
\mathcal B_1^0 = \mathcal B^0 \cap  \mathcal B_1
= \{f \ \text{\hol\ on} \ \ov \D, \ \ f(z) \ne 0 \ \text{on} \
\ov \D; \ \ \|f_0\|_\iy < 1\}.
\end{equation}
By Lemma 3.1, $\mathcal B_1^0$ is a domain in
$H^\iy = H^\iy (\D)$ located
in the unit ball $H_1^\iy$ of this space. Note also that
$$
\sup_{\mathcal B_1^0} |c_n(f)| = \sup_{\mathcal B_1} |c_n(f)|,
$$
and by the maximum principle, each extremal
function $f_0 \in \mathcal B_1$ maximizing $|c_n(f)|$ satisfies
$\|f_0\|_\iy = 1$.

\medskip
\noindent
{\em (b)} \ Since any function $f \in \mathcal B$ does not vanish in $\D$,
the function
$\log f$ is well-defined in a neighborhood of the origin, taking
the principal branch of the logarithmic function, and after \hol \
extension generates a single valued \hol \ function
\be\label{3.5}
\mathbf j_f(z) = \log f(z): \ \D \to \C_{-} := \{w \in \C: \ \Re w < 0\}.
\end{equation}
In fact, we lift $f$ to the universal cover
$$
\C_{-} \to \D_{*} = \D \setminus \{0\}
$$
with the \hol \ universal covering map $\exp$ (cf. Lemma 3.13).

Each such $\mathbf j_f$ satisfies
$$
\sup_\D (1 - |z|^2)^\a |\log f(z)| \le
\sup_\D ((1 - |z|^2)^\a (\log |f(z)| + |\arg f(z)|) < \iy
$$
for any $\a > 0$.
We embed the set $\mathbf j \mathcal B$ into in the Banach space $\B$ of
hyperbolically bounded \hol \ functions on the disk $\D$ with
norm
$$
\|\psi\|_\B = \sup_\D (1 - |z|^2)^2 |\psi(z)|.
$$
This space is dual to the space
$A_1 = A_1(\D)$ of integrable \hol \ functions on $\D$
with $L_1$-norm, and every continuous linear functional
$l$ on $A_1$ can be represented, uniquely, as
\be\label{3.6}
l(\vp) = \langle \psi, \vp\rangle_\D := \iint\limits_\D
(1 - |z|^2)^2 \ov{\psi(z)} \vp(z) dx dy
\end{equation}
with some $\psi \in \B$ (see \cite{Be}).

We want to investigate the geometrical properties of the image
$\mathbf j \mathcal B_1^0$. First of all, we have

\begin{lem} The functions $\mathbf j_f \in \mathbf j \mathcal B$
fill a convex
set in $\B$. Similarly, the subset $\mathbf j \mathcal B_1^0$ is also
convex in $\B$.
\end{lem}

\medskip
\noindent
{\bf Proof}. For any two distinct points
$\psi_1 = \mathbf j f_1, \ \psi_2 = \mathbf j f_2$,
the points of joining interval
$\psi_t = t \psi_1 + (1 - t) \psi_2$ with $0 \le t \le 1$
represent the functions
$\mathbf j f_t = \log (f_1^t f_2^{1-t})$ and the product
$f_1^t(z) f_2^{1-t}(z) \ne 0$ in $\D$ (taking again the
principal branch of logarithm). Hence, this interval lies
entirely in $\mathbf j \mathcal B$.
The proof for $\mathbf j \mathcal B_1^0$ is similar.

\begin{lem} The map $\mathbf j$ is a \hol \ embedding of the domain
 $\mathcal B_1^0$ into the space $\B$ carrying this domain
 onto a \hol \ Banach manifold modelled by $\B$.
\end{lem}

\medskip
\noindent
{\bf Proof}. The map $\mathbf j: \ f \to \log f$ is one-to-one and
continuous on $\mathcal B_1^0$.
To check its complex holomorphy, observe that for any
$f \in \mathcal B_1^0, \ h \in H^\iy(\D)$ and sufficiently
small $|t|$ (letting $\mathbf j(f) = \mathbf j_f$),
$$
\mathbf j(f+th) - \mathbf j(f)
= \log \Bigl( 1 + t \fc{h}{f} \Bigr)
= t \fc{h}{f} + O(t^2),
$$
with uniformly bounded remainder for $\|h\|_\iy \le c < \iy$.
This means that the directional derivative of $\mathbf j$ at $f$
equals $h/f$ and also belongs to $H^\iy(\D)$.

In a similar way, one obtains that the inverse map
$\mathbf j^{-1}: \ \psi \to \e^{\psi}$ is \hol \
on intersections of a neighborhood of $\psi$ in $\B$
with complex lines $\psi + t \om$ in $\mathbf j \mathcal B_1^0$.
The lemma is proved.

Holomorphy in Lemma 3.3 is a special case of general results
on properties of bounded complex Banach functions (see
Lemma 3.12). It implies that both complex structures on
$\mathbf j \mathcal B_1^0$ endowed by norms on $H^\iy$ and on $\B$
are equivalent.

\bk
\noindent
{\bf 2. Complex geometry of sets $\mathbf j \mathcal B_1^0$
and $\mathcal B_1^0$}.

\medskip
As a domain on a complex manifold modelled by $\B$, the set
$\mathbf j \mathcal B_1^0$
admits the invariant \Ko \ and \Ca \ metrics. Our goal is
to show that the geometric features of this set are similar
to bounded convex domains in Banach spaces.

\begin{prop} (i) The \Ko \ and \Ca \ distances on
$\mathbf j \mathcal B_1^0$ and
the corresponding differential metrics are equal:
\be\label{3.7}
\begin {aligned}
d_{\mathbf j \mathcal B_1^0}(\psi_1, \psi_2) &=
c_{\mathbf j \mathcal B_1^0} (\psi_1, \psi_2) =
\inf \{d_\D (h^{-1}(\psi_1), h^{-1}(\psi_2)):
\ h \in \Hol(\D, \mathbf j \mathcal B_1)\}, \\
\mathcal K_{\mathbf j \mathcal B_1^0}(\psi, v) &=
\mathcal C_{\mathbf j \mathcal B_1^0}(\psi, v) \quad
\text{for all} \ \ (\psi, v) \in T(\mathbf j B_1^0).
\end{aligned}
\end{equation}

(ii) Every distinct pair of points $(\psi_1, \psi_2)$
in $\mathbf j \mathcal B_1^0$ can be joined by a complex
$c_{\mathbf j \mathcal B_1^0}$-geodesic.
\end{prop}

\medskip
\noindent
{\bf Proof}. The equality (3.7) follows from the property (ii).
We establish this property in two steps.

$(a)$ \ First take the $\epsilon$-blowing up of
$\mathbf j \mathcal B_1^0$, that is, we consider the sets
$$
U_\epsilon = \bigcup_{\psi \in \mathbf j \mathcal B_1^0} \
\{\om \in \B: \ \|\om - \psi \|_\B < \epsilon\},
\quad \epsilon > 0.
$$

For these sets, we have

\begin{lem} Every set $U_\epsilon$ is a (bounded) convex domain
in $\B$, and its weak$^*$-closure in $\sigma(\B, A_1)$ is
compact.
\end{lem}

\medskip
\noindent
{\bf Proof}.
The openness and connectivity of $U_\epsilon$ are trivial.
Let us check convexity.
Take any two distinct points $\om_1, \om_2$ in $U_\epsilon$
and consider the line interval
\be\label{3.8}
\om_t = t \om_1 + (1 - t) \om_2, \quad 0 \le t \le 1,
\end{equation}
joining these points. Since, by definition of $U_\epsilon$,
each point  $\om_n \ (n = 1, 2)$ lies in the ball
$B(\psi_n, \epsilon)$ centered at $\psi_n$ with radius
$\epsilon$, and the interval
$\{\psi_t = t \psi_1 + (1 - t) \psi_2\}$
lies in $\mathbf j B_1^0$, we have, for all $0 \le t \le 1$,
$$
\om_t - \psi_t = t(\om_1 - \psi_1) + (1 - t) (\om_2 - \psi_2)
$$
and
$$
\|\om_t - \psi_t\| \le t\|\om_1 - \psi_1\| +
(1 - t) \|\om_2 - \psi_2\| < \epsilon,
$$
which shows that the interval (3.8) lies entirely in $U_\epsilon$.

To establish $\sigma(\B, A_1)$-compactness of
the closure $\ov U_\epsilon$,
note that weak$^*$ convergence of the functions $\om_n \in \B$
to $\om$ implies the uniform convergence of these functions on
compact subsets of $\D$.
It suffices to show that for any bounded sequence
$\{\om_n\} \subset \B$ we have the equality
\be\label{3.9}
\lim\limits_{n\to \iy} \iint\limits_\D
\fc{(1 - |\z|^2)^2 \om_n(\z)}{\z - z} d \xi d \eta =
\iint\limits_\D
\fc{(1 - |\z|^2)^2 \om(\z)}{\z - z} d \xi d \eta, \quad z \in \D^*,
\end{equation}
because the functions $w_z(\z) = 1/(\z - z)$ span a dense subset of $A_1(\D)$.
But if
$$
\sup_\D (1 - |\z|^2)^2 |\om(\z)| < M < \iy \quad
\text{for all} \ \ n,
$$
the equality(3.9) follows from Lebesgue's theorem on dominant
convergence. Lemma follows.

\bk
$(b)$ \ We proceed to the proof of Proposition 3.4 and first
establish the existence of complex geodesics in domains
$U_\epsilon, \ \epsilon < \epsilon_0$.
Convexity of these domains allows us to use the arguments applied in
\cite{Di} in the proof of Proposition 2.3.

Let $\om_1$ and $\om_2$ be distinct points in $U_\epsilon$.
By Proposition 2.1,
$$
d_D(\om_1, \om_2) = c_D (\om_1, \om_2) =
\inf \{d_\D (h^{-1}(\om_1), h^{-1}(\om_2)): \ h \in
\Hol(\D, U_\epsilon)\};
$$
hence there exists the sequences
$\{h_n\} \subset \Hol(\D, U_\epsilon)$  and
$\{r_n\}, \ 0 < r_n < 1$, such that
$h_n(0) = \om_1$ and $h_n(r_n) = \om_2$ for all $n, \
\lim\limits_{n\to \iy} r_n = r <1$ and
$c_{U_\epsilon} (\om_1, \om_2) = d_\D(0, r)$.
Let $h_n(t) = \sum\limits_{m=0}^\iy a_{m,n} t^m$
for all $t \in \D$ and $n$.

Take a disk $\D_R =\{|z| < R\}$ containing $U_\epsilon$. The
Cauchy inequalities imply
$\|a_{n,m}\|_\B \le R$ for all $n$ and $m$.
Passing, if needed, to a subsequence of $\{h_n\}$, one can suppose that
for a fixed $m$, the sequence $a_{n,m}$ is weakly$^*$ convergent
to $a_m \in \B$ as $n \to \iy$, that is
$$
\lim\limits_{n\to \iy} \langle a_{n,m}, \vp\rangle_\D =
\langle a_m, \vp\rangle_\D \quad \text{for any} \ \ \vp \in A_1.
$$
Hence $h(t) = \sum\limits_{m=0}^\iy a_m t^m$
defines a \hol \ function from $\D$ into $\B$.
Since $a_{n,0} = \z$ for all $n$, we have $h(0) = \om_1$.

Now, let $\a, \ 0 <\a < 1$, and $\ve > 0$ be given. Choose $m_0$
so that
$$
r \sum\limits_{m=m_0}^\iy \a^m < \ve.
$$
If $\vp \in A_1, \ \|\vp\| = 1$, then
$$
\sup_{|t|\le \a} | \langle h_n(t) - h(t), \vp \rangle_\D|
\le \sum\limits_{m=1}^{m_0-1}
|\langle a_{n,m} - a_m, \vp \rangle_\D|
+ 2 r \sum\limits_{m=m_0}^\iy \a^m
$$
for all $n$, which implies that $h_n$ is convergent to $h$
in $\sigma (\B, A_1)$ uniformly on compact subsets of $\D$
as $n \to \iy$. Since $\ov D$ is $\sigma (\B, A_1)$ compact,
$h(\D) \subset \ov D$, and since $h(0)  \in D$, it follows that
$h(\D) \subset D$. For $r < r^\prime < 1$,
\be\label{3.10}
\om_2 = h_n(r_n) = \fc{1}{2 \pi i} \int\limits_{|t| = r^\prime}
\fc{h_n(t)dt}{t - r_n} \to \fc{1}{2 \pi i} \int\limits_{|t| = r^\prime}
\fc{h(t)dt}{t - r} = h(r)
\end{equation}
as $n \to \iy$. Hence,
$$
d_T(0, r) = c_D (\om_1, \om_2) = c_D(h(0), h(r)),
$$
and $h$ is simultaneously complex $c_{U_\epsilon}$ and
$d_{U_\epsilon}$ geodesics.

There exists a \hol \ map $g: \D \to U_\epsilon$ such that
for any two points $t_1, t_2 \in \D$,
\be\label{3.11}
d_\D(t_1, t_2) = d_{U_\epsilon}(g(t_1), g(t_2)) =
c_{U_\epsilon}(g(t_1), g(t_2)),
\end{equation}
and for any pair $(t, v), \ t \in \D, \ v \in \C$,
\be\label{3.12}
\mathcal K_{U_\epsilon}(g(t), d g(t) v) = \fc{|v|}{1 - |t|^2}.
\end{equation}

\bk
$(c)$ \ Let now $\om_1$ and $\om_2$ be two distinct points in
$\mathbf j \mathcal B_1^0$.
Choose a decreasing sequence $\{\epsilon_n\}$ approaching
zero and take for every $n$ a complex geodesic $h_n = h_{U_{\epsilon_n}}$
joining these points in $U_{\epsilon_n}$, which was constructed in the
previous step. Let $g_n = g_{U_{\epsilon_n}}$ be the corresponding map
$\D \to U_{\epsilon_n}$ which provides the equalities (3.11), (3.12).
Since $d_\D$ is conformally invariant, one can take $g_n$ satisfying
$g_n^{-1}(\om_1) = 0, \ g_n^{-1}(\om_2) = r_n \in (0, 1)$. Then the
inequalities
$$
d_{U_{\epsilon_n}}(\om_1, \om_2) \le d_{U_{\epsilon_m}} (\om_1, \om_2)
\le d_{\mathbf j B_1} (\om_1, \om_2) \quad \text{for} \ \ m > n
$$
imply $r_n \le r_m \le r_{*} < 1$, where
$d_\D(0, r_{*}) = d_{\mathbf j \mathcal B_1}(\om_1, \om_2)$.
Hence, there exists
$\lim\limits_{n\to \iy} r_n = r^\prime \le r_{*}$.

The sequence $\{g_n\}$ is $\sigma(\B, A_1)$-compact and similar
to (3.10) the weak$^*$ limit of $g_n$ is a function
$g \in \Hol(\D, \mathbf j \mathcal B_1^0)$
which determines a complex geodesic for both \Ko \ and \Ca \ distances on
${\mathbf j \mathcal B_1^0}$ joining the points $\om_1$
and $\om_2$ inside this set.
Proposition 3.4 is proved.

An important consequence of Proposition 3.4 is that the initial
domain $\mathcal B_1^0$ in $H^\iy$ has similar complex geometric
properties, since the embedding $\mathbf j$ is biholomorphic.
We present it as

\begin{prop} (i) The \Ko \ and \Ca \ distances on domain
$\mathcal B_1^0$ and
the corresponding differential metrics are equal:
\be\label{3.13}
\begin {aligned}
d_{\mathcal B_1^0}(f_1, f_2) &=
c_{\mathcal B_1^0} (f_1, f_2) =
\inf \{d_\D (h^{-1}(f_1), h^{-1}(f_2)):
\ h \in \Hol(\D, \mathcal B_1)\}, \\
\mathcal K_{\mathcal B_1^0}(f, v) &=
\mathcal C_{\mathcal B_1^0}(f, v) \quad
\text{for all} \ \ (f, v) \in T(B_1^0).
\end{aligned}
\end{equation}

(ii) Every two points $f_1, f_2$ in $\mathcal B_1^0$ can be joined
by a complex geodesic.
\end{prop}

\bk
\noindent
{\bf 3. Finsler metric generated by functional $c_n(f)$}.

\medskip
We proceed to the proof of Theorem 1.1.
It will be convenient to regard the free coefficients
$c_0(f)$ also as elements of $\mathcal B_1^0$, which are
constant on $\D$. Note that $0 < |c_0(f)| < 1$. Denote
$$
\sup_{f\in \mathcal B_1^0} |c_n(f)| =
\sup_{f\in \mathcal B_1} |c_n(f)| = M_n
\quad (M_n \le 1),
$$
and consider, for a fix integer $n > 1$, the functional
\be\label{3.14}
J(f) = \Big\vert \fc{c_n(f)}{M_n} \Big\vert^{1/n}: \
\mathcal B_1^0 \to [0, 1).
\end{equation}
It is logarithmically \psh \ on $\mathcal B_1^0$,
taking the values on $[- \iy, 0)$, with
$\log J(f) \to 0$ as $f$ tends to the boundary of
$\mathcal B_1^0$.

Our goal is to show that $\log J$ is dominated on
$\mathcal B_1^0$ by the pluricomplex Green function of this domain,
namely,
\be\label{3.15}
\log J(f) \le g_{\mathcal B_1^0}(c_0(f), f).
\end{equation}
This will be established in several steps.
The proof is geometric and involves the differential metrics.
We construct on each \hol \ disk in $\mathcal B_1^0$
a \sh \ Finsler metric naturally generated by $J$ and compare
this metric with the canonical \Ko \ metric.

The estimate (3.15) trivially holds for the points of the zero-set
of our functionals
\be\label{3.16}
Z_J = \{f \in \mathcal B_1^0: \ J(f) = 0\}.
\end{equation}
Note that this set contains the disk filled in $\mathcal B_1^0$
by the constant functions $c_0(f)$,  which can be identified with
the punctured disk $\D \setminus \{0\}$.

The uniqueness theorem for \hol \ functions (in Banach spaces)
implies that this zero-set is nonwhere dense on $\mathcal B_1^0$
(in the sense that its complement
$\mathcal B_1^0 \setminus Z_J$ is open and dense everywhere).
This follows also from a theorem from \cite{Kr1} on the existence
of special \qc \ deformaions $\om$ of the plane, which are conformal
on a given set $E$ of positive two-dimensional Lebesgue's measure
and take, with their derivatives, the prescribed values,
see [Kr1, Ch. 4]. One can compose any function
$f \in \mathcal B_1^0$ with appropriate deformations of such
kind, that are conformal, for example, in the complement of
a disk $\{|z - z_0| < r$ located sufficiently far from the origin,
and get the composite maps $\om \circ f$ whose coefficients
$c_n(\om \circ f)$ range in a whole neighborhood $0 < |w| < \ve$.
Then $\om \circ f$ provide the points from
$\mathcal B_1^0 \setminus Z_J$.

\bk
Consider first the \hol \ disks $h: \D \to \mathcal B_1^0$ with
nonconstant \hol \ $h$ which touch the zero-set (3.16) only at one
point. We call such disks {\bf distinguished}. One can assume
that this common point is $h(0)$.

Let $\mathcal D = h(\D)$ be such a disk.
Take the restriction of $J(f)$ to $\mathcal D$ and
consider its root
$$
g(\z) = [J \circ h(\z)]^{1/n};
$$
this root is an $n$-valued function on $\mathcal D$,
with a single algebraic branch point at $\z = 0$. Take a
single-valued branch of this function in a neighborhood
$U_0 \subset \mathcal D$ of a point $\z_0 \ne 0$ and apply
the selected branch to pulling
back the hyperbolic metric $\ld_{\hyp}$ on $\D$ to this neighborhood
$U_0$.
Extending this branch analytically, one produces
a conformal metric $ds = \ld_J(\z) |d\z|$ on the whole disk
$\mathcal D$, with
\be\label{3.17}
\ld_J(\z) = g^* \ld_{\hyp}(\z) = \fc{|g^\prime(\z)|}{1 - |g(\z)|^2}.
\end{equation}
This metric does not depend on the choices of an initial branch and
of $U_0$. It is logarithmically \sh \ on $\mathcal D$,  and its
Gaussian curvature $\kappa_{\ld_g}$ equals $- 4$ at noncritical
points of the extension of $g$. This provides that the
curvature is less than or equal $- 4$ on $\mathcal D$
in both supporting and potential senses and as generalized curvature
via (2.8).

\bk
\begin{lem} On any $d_{\mathcal B_1^0}$-geodesic disk $\mathcal D$,
the metric (3.17) is dominated by the differential \Ko \ metric
$\ld_{\mathcal K_{\mathcal B_1^0}}$,
\be\label{3.18}
\ld_J(\z) \le \ld_{\mathcal K_{\mathcal B_1^0}}(\z),
\end{equation}
and if the equality holds here for one value of $\z \ne 0$, then it holds
identically.
\end{lem}

\medskip
\noindent
{\bf Proof}. Consider first the distinguished geodesic disks
$\mathcal D = h(\D)$.
On such disks, the differential \Ko \
metric $\ld_{\mathcal K_{\mathcal B_1^0}}$ is equal to
hyperbolic metric of the unit disk,
Since the curvature of $\ld_J$ is at most $- 4$ at noncritical
points of $\mathcal D$ in the supporting sense, the inequality (3.18)
follows from the classical Ahlfors-Schwarz lemma (see \cite{Ah},
\cite{He}, \cite{Mi}, \cite{Ro}).

An arbitrary geodesic disk in  $\mathcal B_1^0$ can be strongly
(in the norm of $H^\iy$) approximated by distinguished disks, and
such approximation preserves the inequality (3.18). This completes
the proof.

\bk
We must now pass from the inequality (3.18) for infinitesimal
metrics to global distances, which requires the reconstruction
of the initial functional $J(f)$ from the generated metric.
This is rather simple for distinguished geodesic disks.

\begin{lem} On any distinguished geodesic disk
$h: \ \D \to \mathcal B_1^0$, we have for each $r < 1$
the equality
\be\label{3.19}
\tanh [J(f)] =
\tanh^{-1}[\dl(J \circ h(0), J \circ h(r))]
= \int\limits_0^r \ld_J \circ h(t) dt,
\end{equation}
where
$$
\dl (\z_1, \z_2) = (\z_2 - \z_1)/(1 - \ov \z_1 \z_2).
$$
\end{lem}

\medskip
\noindent
{\bf Proof}. Since any geodesic disk is holomorhically
isometric to the
hyperbolic plane modelled by $\D$, one can write
\be\label{3.20}
\tanh^{-1}[\dl(J \circ h(0), J \circ h(r))] =
\int\limits_{J \circ h(0)}^{J \circ h(r)}
\fc{|d t|}{1 - |t|^2}
= \int\limits_0^r \ld_{J \circ h}(t) |d t|
\quad (0 < r < 1),
\end{equation}
which is equivalent to (3.19).
Indeed, one can subdivide the hyperbolic interval
$[\dl(J \circ h(0), J \circ h(r))]$
into subintervals, taking a finite partition
$$
c_0 <r_1 < \dots < r_{m-1} < r_m = J \circ h(r)
$$
so that on each $[r_{s-1}, r_s]$ the map $J \circ h$
is injective, and apply to these subintervals the equalities
similar to (3.20).

\bk
Note also that {\em if a geodesic disk $h(\D)$ is not distinguished,
but does not lie
entirely in $Z_J$, then the equality (3.19) holds for
a sufficiently small $r < 1$}, for which the initial equality
(3.20) remains valid. The same holds, in view of (2.6), for
any compact subset
of the disk $h(\D)$ which does not contact the zero set $Z_J$.

\bk
Lemmas 3.7 and 3.8, together with Proposition 3.6, imply the desired
estimate (3.15) which controls the behavior of $J$ on $\mathcal B_1^0$.
In view of its importance, we present this as a separate lemma.

\begin{lem} The inequality (3.15) holds at every point
$f \in \mathcal B_1^0$.
\end{lem}

\medskip
\noindent
{\bf Proof}. The case $J(f) = 0$ is trivial, so
we have to establish the inequality (3.15) only for the points
$f$ with $J(f) \ne 0$.

Lemmas 3.7 and 3.8 imply that the growth of $J$ on the distinguished
geodesic disks is estimated by
$$
J(f) = O(d_{\mathcal B_1^0} (c_0(f), f))
= O(\|f - c_0(f)\|_{H^\iy}),
$$
uniformly on compact subsets of these disks. This estimate
provides that $\log J(f)$ is an admissible \psh \ function
for comparison with Green's function
$g_{\mathcal B_1^0}(c_0(f), f)$.
Now the maximality of $g_{\mathcal B_1^0}(c_0(f), f)$
among \psh \ functions
with logarithmic growth near the pole $c_0$ implies that $\log J(f)$
is dominated by Green's function $g_{\mathcal B_1^0}(c_0(f), f)$.

Further, for any given function $f \in \mathcal B_1^0$,
the point $c_0(f)$ belongs to the zero-set $Z_J$.
Using approximation in $\mathcal B_1^0$, similar to above,
one can extend the inequality (3.15) to all complex geodesic disks in
$\mathcal B_1^0$ which touch $Z_J$ at this point.
By continuity, the functional $\log J$ is \sh \ on every
such disk, while the relations (2.6), (3.18) and (3.19) preserve
the required logarithmic order of the growth of $J$ near
its zero set. Lemma follows.

\bk
\noindent
{\bf 4. \ Homotopy}.

\medskip
For any $f \in \mathcal B_1^0$, one can define complex \hol \
homotopy
\be\label{3.21}
f_t(z) = f(t z) = c_0 + c_1 t z + \dots: \
\D \times \D \to \D_{*}
\end{equation}
connecting $f$ with $c_0(f)$ in $\mathcal B_1^0$.
Due to (3.14), our functional $J$ is homogeneous
with respect to this isotopy with degree $1$ in the
following sense:
\be\label{3.22}
J(f_t) = |t| J(f).
\end{equation}

We shall need also the following simple fact concerning the
homotopy functions.

\begin{lem} The pointwise map $t \mapsto f_t$
given by (3.23) determines a \hol \ map $\chi_f: \D \to H^\iy$.
\end{lem}

This lemma is a rather special case of bounded \hol \ functions
in Banach spaces with sup norm, given by the following lemma (cf. \cite{Ea},
\cite{Ha}, \cite{Kr2}).

\begin{lem} Let $E, \ T$ be open subsets of complex Banach spaces
$X, Y$ and $B(E)$ be a Banach space of \hol \ functions
on $E$ with sup norm. If $\vp(x, t)$ is a bounded map $E \times T \to B(E)$
such that $t \mapsto \vp(x, t)$ is \hol \ for each $x \in E$,
then the map $\vp$ is \hol.
\end{lem}

We briefly outline the proof.
Holomorphy of $\vp(x, t)$ in $t$ for fix $x$ implies the existence of
complex directional derivatives
$$
\vp_t^\prime(x,t) = \lim\limits_{\z\to 0}
\fc{\vp(x, t + \z v) - \vp(x, t)}{\z} =
\fc{1}{2 \pi i} \int\limits_{|\xi|=1} \fc{\vp(x, t + \xi v)}{\xi^2} d \xi.
$$
On the other hand, the boundedness of $\vp$ in sup norm
provides the uniform estimate
$$
\|\vp(x, t + c \z v) - \vp(x, t) - \vp_t^\prime(x,t) c v\|_{B(E)}
\le M |c|^2,
$$
for sufficiently small $|c|$ and $\|v\|_Y$.

\bk
\noindent
{\bf 5. \ Covering maps}.

\medskip
For each $f_t$, there exists, by Proposition 3.4, a
complex geodesic in $\mathcal B_1^0$ joining $f_t$ with
$c_0(f)$. It is a \hol \ geodesic disk isometric to the
hyperbolic plane $\mathbb H^2$, i.e., with
the same hyperbolic geometry as on $\mathbb H^2 = \D$.
We need to estimate quantitatively the behavior of
the distance  $d_{\mathcal B_1^0} ({f_t}, c_0)$
when $t \to 0$.

\begin{lem} Every function $f \in \mathcal B_1^0$ admits
factorization
\be\label{3.23}
f(z) = \kappa_0 \circ \wh f(z),
\end{equation}
where $\wh f$ is a \hol \ map of the disk $\D$ into itself
(hence, from $H_1^\iy$) and $\kappa_0$ is the function (1.2).
\end{lem}

\medskip
\noindent
{\bf Proof}. Due to a general topological theorem,
any map $f: M \to N$, where $M, N$ are manifolds, can be
lifted to a covering manifold $\wh N$ of $N$,
under appropriate relation between the fundamental group
$\pi_1(M)$ and a normal subgroup of $\pi_1(N)$ defining
the covering $\wh N$ (see, e.g, [Ma]).
This construction produces a map $\wh f: M \to \wh N$
satisfying
\be\label{3.24}
f = p \circ \wh f,
\end{equation}
where $p$ is a projection $\wh N \to N$. The map
$\wh f$ is determined
up to composition with the covering transformations
of $\wh N$ over $N$.
For \hol \ maps and manifolds the lifted map is also
\hol.

In our special case, $\kappa_0$ is a \hol \ universal
covering map $\D \to \D_{*} = \D \setminus \{0\}$, and
the representation (3.24) provides the equality (3.23)
with the corresponding $\wh f$ determined up to
covering transformations of the unit disk compatible
with the covering map $\kappa_0$.

\bk
This lemma relates to Lemma 3.2.
As a simple corollary of Lemma 3.12 one obtains

\begin{lem} For any $f \in \mathcal B_1^0$,
\be\label{3.25}
|c_1| \le 2/e,
\end{equation}
with equality only for the rotations
$e^{i \a_1} \kappa_0(e^{i \a} z)$ of $\kappa_0$.
\end{lem}

\bk
As was mentioned in the introduction, this bound is known.
We reprove it and will use also later the arguments applied
in the proof.

\medskip
\noindent
{\bf Proof}. One only needs to show that (3.25) holds for each
composition of $\kappa_0$ with the M\"{o}bius (fractional linear)
automorphisms $\g$ of
of the unit disk $\D$, i.e., that $\kappa_0$ (and any its rotation)
maximizes $|c_1|$ among the \hol \ universal  covering maps
$\D \to \D_{*}$.
Then, for any $f \in \mathcal B_1^0$, taking $\g$ with
$\g(0) = \wh f(0) = a, \ \g(1) = \wh f(1)$, where
$\wh f$ is the covering map of $f$ in (3.23), one obtains from
Schwarz' lemma,
$$
|c_1(f)| = |\kappa_0(a)| |\wh f^\prime(0)|
\le |c_1(k_\g(0))| = |\kappa_0(a)| |\g^\prime(0)| \le \fc{2}{e}.
$$

Using the rotations about the origin $z = 0$, we can restrict
ourselves by $\g$ whose compositions with
$\sigma(z) = (z - 1)/(z + 1)$ assumes the form
$$
\sigma \circ \g (z) = e^{i b} \fc{z - 1}{z + e^{i a}}
\quad \text{with} \ \ a, b \in [0, 2 \pi].
$$
We regard here the disk $\D$ as a lune
with vertices at the points $z = 1$ and $z = e^{i a}$
and with opening angles equal to $\pi$.
Then
$$
(\kappa_0 \circ \g)^\prime(z) = e^{\sigma \circ \g(z)}
\fc{A}{(z +  e^{i a})^2}, \quad  A = e^{i b} (1 + e^{i a}),
$$
and
$c_1(\kappa_0 \circ \g)| = |(\kappa_0 \circ \g)^\prime(0)| = |A|$,
with equality only for $a = 0$. Lemma follows.

We will denote the M\"{o}bius group of $\D$ by $\Mob (\D)$
and put
$$
\g^* \kappa_0 := \kappa_0 \circ \g.
$$

\noindent
{\bf 5. \ Estimates for the \Ko \ distance on $\mathcal B_1^0$}.

\medskip
\begin{prop} For any $f =  \in \mathcal B_1^0$,
we have the equality
\be\label{3.26}
d_{\mathcal B_1^0} (f, c_0)
= \inf \{d_{H_1^\iy} (\wh f, \wh c_0): \ \kappa_0 \circ \wh f = f\};
\end{equation}
moreover, there exists a map
$\wh f^*(z) = c_0^* + c_1^* z + \dots$ covering $f$, on which
the infimum in (3.26) is attained, i.e.,
\be\label{3.27}
d_{\mathcal B_1^0} (f, c_0) = d_{H_1^\iy} (\wh f^*, \wh c_0^*).
\end{equation}
\end{prop}

\medskip
\noindent
{\bf Proof}.
It is well-known (and rather simple) that a
complex geodesic in the unit ball $B(0, 1)_{\mathbf X}$ of
a complex Banach space $\mathbf X$, joining its center
$\mathbf 0$ with a point $\x \ne 0$, is a \hol \ isometry
$\D \to B(0, 1)_{\mathbf X}$ determined by the map
$$
\z \mapsto \z \x/\|\x\|
$$
and that the \Ko \ and \Ca \ distances in $B(0, 1)_{\mathbf X}$
between these points are equal to
\be\label{3.28}
d_{B(0, 1)_{\mathbf X}}(\mathbf 0, \x) = d_\D (0, \|\x\|)
= \tanh^{-1} \|\x \|.
\end{equation}
We apply this to functions
$\wh f(z) = \wh c_0 + \wh c_1 z + \dots$ from $H_1^\iy$.
The corresponding functions
$$
\wh g_f(z) := \fc{\wh f(z) - \wh c_0}{1 - \ov{\wh c_0} \wh f(z)},
$$
also belong to $H_1^\iy$, and by (3.28),
\be\label{3.29}
d_{H_1^\iy}(\wh g_f, \mathbf 0) = d_\D (\|\wh g_f\|_\iy, 0)
= \tanh^{-1} (\|(\wh f - \wh c_0)/(1 - \ov{\wh c_0} \wh f)\|_\iy).
\end{equation}
Since for a fixed $\wh c_0$ (regarded again as a constant
function on $\D$), the map
$$
\wh g \mapsto \fc{\wh g + \wh c_0}{1 + \ov{(c_0)} \wh g}
$$
with $\wh g$ running over the ball $H_1^\iy$ is a biholomorphic
isometry of this ball, the map
\be\label{3.30}
\om (\z) = \fc{\z \wh g_f/\|\wh g_f\|_\iy + \wh c_0}{1 +
\ov{\wh c_0}  \z \wh g_f/\|\wh g_f\|_\iy}: \ \D \to H_1^\iy
\end{equation}
carries out the complex geodesic
$\z \mapsto \z \wh g_f/\|\wh g_f\|_\iy$
into a complex geodesic in $H_1^\iy$ passing through
the points $\wh c_0$ and $\wh f$. The point $\wh f$ is obtained
by (3.30) on the value $\z = \|\wh g_f\|_\iy$.

Now observe that the universal covering map
$\kappa_0: \D \to \D_{*}$
extends by the equality (3.23) to all $\wh f \in H_1^\iy$ and
this extension
induces a \hol \ map of the unit ball $H_1^\iy$ into
domain $\mathcal B_1^0$. Such maps cannot expand the invariant
distances; thus
$$
d_{\mathcal B_1^0}(f, c_0) = d_{\mathcal B_1^0}
(\kappa_0 \circ \wh f, \kappa_0(\wh c_0))
\le d_{H_1^\iy}(\wh f, \wh c_0),
$$
and
\be\label{3.31}
d_{\mathcal B_1^0}(f, c_0)
\le \inf_{\wh f} d_{H_1^\iy}(\wh f, \wh c_0),
\end{equation}
where the infimum is taken over all covers $\wh f$ of $f$.

Our goal is to show that in fact one has the equality in (3.31)
and that the infimum in (3.31) is attained.
We assume, the contrary, i.e., that
$$
d_{\mathcal B_1^0}(f, c_0) < \inf_{\wh f} d_{H_1^\iy}(\wh f, \wh c_0),
$$
and apply Proposition 3.6.
This proposition provides the existence of a complex geodesic
$h: \ \D \to \mathcal B_1^0$ joining there the points
$c_0$ and $f$, and it follows from the above relations that
\be\label{3.32}
d_{\mathcal B_1^0}(f, c_0) = d_\D(\z_1, \z_2)
< \tanh^{-1} \|(f - c_0)/(1 - \ov c_0 f)\|_\iy,
\end{equation}
where $\z_1 = h^{-1}(c_0), \ \z_2 = h^{-1}(f)$.
Lifting the map $h$ by (3.23) to a \hol \ map $\wh h$ of the unit disk
into itself, one gets the points
$\wh h(\z_1), \wh h(\z_2)$ in $\D$, which lie in the fibers
over $c_0$ and $f$, respectively, while the relations (3.32) imply
$$
d_{H_1^\iy}(\wh h(\z_1), \wh h(\z_2)) = d_\D(\z_1, \z_2)
< \tanh^{-1} \|(f - c_0)/(1 - \ov c_0 f)\|_\iy.
$$
This inequality contradicts (3.29), which completes the
proof of Proposition 3.14.

\medskip
This proposition provides explicitly some complex geodesics
in $\mathcal B_1^0$.
Note also that the arguments in the proof above remain in
force by
replacing there $\kappa_0$ by a universal covering map
$\g^* \kappa_0: \ \D \to \D_{*}$ with a fixed
$\g \in \Mob (\D)$.

\bk
We now choose the factor $\wh f$ in (3.23) so that $\wh f(0) = 0$
and fix such $\wh f$. Accordingly, we must replace $\kappa_0$ by
$\g^* \kappa_0 = \kappa_0 \circ \g$ with $\g \in \Mob(\D)$,
which satisfies
$$
\g(0) = \kappa_0^{-1}(c_0),
$$
and taking the point $\kappa_0^{-1}(c_0)$ in the closure of the
fundamental triangle of the Fuchsian group $\G(\D, \D_{*})$
which uniformizes the punctured disk $\D_{*}$ in $D$
(that means $\D_{*}$ is represented as factor
$D/\G(\D, \D_{*})$ up to conformal equivalence;
the desired conformal map is produced by $\kappa_0$).
Then the representation (3.23) assumes the form
\be\label{3.33}
f(z) = (\g^* \kappa_0) \circ \wh f(z).
\end{equation}
Note that $\g$ is determined up to rotations about the origin,
which is not essential for dilatation.

\bk
As a corollary of Proposition 3.14, one obtains the following
result, which we precede by some remarks.
Denote again the coefficients of the covering maps $\wh f$ in (3.33)
by $\wh c_n$ and note that, in view of the assumption $\wh c_0 = 0$,
$$
\wh f(z) = \wh c_1 z + \wh c_2 z^2 + \dots \ .
$$
Accordingly,
$\wh f_t(z) = \wh c_1 z + \wh c_2 z^2 + \dots$
will denote the covers of homotopies (3.21)
for original functions $f \in \mathcal B_1^0$.

It follows from Lemma 3.11 that
$$
\|\wh f_t\|_\iy = |\wh c_1| |t| + O(t^2)
\quad \text{as} \ \ t \to 0,
$$
where the estimate of remainder is in $H^\iy$-norm (thus
uniform for all $|z| < 1$).
If $c_1 = c_2 = \dots = c_{m-1} = 0$
(equivalently, for
$\wh c_1 = \wh c_2 = \dots = \wh c_{m-1} = 0$),
we have
\be\label{3.34}
d_{\mathcal B_1^0} (f_t, c_0) \le d_{H_1^\iy} (\wh f_t, 0)
= |\wh c_m| |t|^m + O(|t|^{m+1});
\end{equation}
here $\wh c_m$ is the first nonvanishing coefficient.

A consequence of Proposition 3.14 mentioned above is the following

\begin{lem} For any function
$f(z) = c_0 + \sum_m^\iy c_n z^n \in \mathcal B_1^0, \ m \ge 1$,
with $c_m \ne 0$ (and $c_0 \ne 0$) we have the sharp asymptotic
estimate
\be\label{3.35}
\begin{aligned}
d_{\mathcal B_1^0} (f_t, c_0) &=\inf_f d_{H_1^\iy} (\wh f_t, 0)
= \inf \{|\wh c_m(\wh f)|: \ \g^* \kappa_0 \circ \wh f = f\}
|t|^m + O(|t|^{m+1})  \\
&= \inf \fc{|c_m|}{(\g^* \kappa_0)^\prime(0)}|t|^m + O(|t|^{m+1}),
\quad t \to 0,
\end{aligned}
\end{equation}
where the infima are taken over all covering maps $\wh f$ of $f$
fixing the origin and all $\g \in \G(\D, \D_{*})$, and
these infima are attained on some pair $(\wh f, \g)$.
\end{lem}

\medskip
In paricular, $\wh \kappa_0(z) = z$, and for this map
the equalities (3.35) result in
$$
d_{\mathcal B_1^0} (\kappa_{0,t}, c_0) = |t| + O(|t|^2),
\quad t \to 0.
$$
This equality shows that the {\em \hol \ disk $\D(\kappa_0)$
filled by the homotopy functions
$\kappa_{0t}(z) = \kappa_0(tz), \ t \in \D$,
is a complex geodesic in} $ \mathcal B_1^0$.

Accordingly, for $\kappa_m(z) = \kappa_0(z^m)$, we have
\be\label{3.36}
d_{\mathcal B_1^0} (\kappa_{m,t}, c_0) = |t|^m + O(|t|^{m+1}),
\quad t \to 0.
\end{equation}
In fact, the remainder terms in the last two equalities
can be omitted.

Note also that the covering map $\wh f$ is a rotation of $\D$ about
the origin only for $f = \kappa_0$, and for all rotations we have
similar results.

\bk
\noindent
{\bf 6. Finishing the proof of Theorem 1.1}.

\medskip
We can now complete the proof of the main theorem.
Let
$$
f^0(z) = c_0^0 + c_1^0 z + c_2^0 z^2 + \dots
$$
be an extremal function maximizing $|c_n| \ (n > 1)$ on
$\mathcal B_1$. Then $|c_n^0| = M_n$ and $J(f^0) = 1$,
and for its homotopy functions
$f_r^0(z) = f^0(r z)$ with $0 < r <1$, which lie in
$\mathcal B_1^0$, we have
$J(f_r^0) = r$.

First we show that any extremal function $f^0$ must satisfy
\be\label{3.37}
 c_1^0 = 0.
\end{equation}
Indeed, assuming $c_1^0 \ne 0$, one derives from Lemma 3.15
the equalities
\be\label{3.38}
d_{\mathcal B_1^0} (f_r^0, c_0^0)
= |\wh c_1^0| r + O(r^2)
= \fc{|c_1^0| r}{|(\g_0^* \kappa_0)^\prime(0)|} + O(r^2),
\quad r \to 0,
\end{equation}
where $\wh c_1^0$ is the first coefficient of a factorizing
function $\wh f^0$ for $f^0$ by (3.33) and $\g$ is the
appropriate M\"{o}bius automorphism of $\D$, on which the
infima in (3.35) are attained.
Combining these equalities with the relations(2.2), (2.6),
connecting the \Ko \ distance and the Green
function $g_{\mathcal B_1^0}(\mathbf 0, f_r^0)$,
and with Lemma 3.9, one obtains
$$
r J(f^0) = r \le  |\wh c_1^0| r + O(r^2),
$$
and therefore,
$$
1 \le |\wh c_1^0| = \fc{|c_1^0|}{|(\g^* \kappa_0)^\prime(0)|}.
$$
By Schwarz's Lemma and Lemma 3.13, such relations can hold only
in the case when $|\wh c_1^0| = 1$, i.e., for $f^0 = \kappa_0$
(up to rotations);
in addition, we must have the equalities
$$
|c_n^0| = |c_1^0| = |c_1(\g_0^* \kappa_0)|= 2/e.
$$
But this is impossible, because in view of Parseval's equality
$\sum_0^\iy |c_n|^2 = 1$ for the boundary function
$\kappa_0(e^{i \theta}), \ \theta \in [0, 2 \pi]$;
in fact, the strict inequality
$$
|c_n(\kappa_0)| < 2/e
$$
holds for any $n > 1$.
This contradiction proves the equality (3.37).

Therefore, the extremal functions of $J$ must be of the form
$f^0(z) = c_0^0 + c_2^0 z^2 + \dots$; equivalently,
\be\label{3.39}
f_t^0(z) = c_0^0 + c_2^0 t^2 z^2 + O(t^3), \quad t \to 0.
\end{equation}
Now the proof of the theorem is continued successively for
$n = 2, 3, \dots$ .

By maximization of the second coefficient $c_2$, the expansion
(3.39) requires to deal with the square functional
$$
J_2(f) = J(f)^2 = |c_2(f)|/M_2
$$
to have homogeneity of degree $2$.
Its comparison with (3.36) (for $m = 2$) provides
$$
J_2(f_r^0) = r^2 \le
\fc{|c_2^0| r^2}{(\g^* \kappa_0)^\prime(0)} + O(r^3),
$$
which implies, similar to (3.38), the equalities
$$
|c_2^0| = \kappa_0^\prime(0) = 2/e.
$$
These equalities yield
$$
f^0(z) = \kappa_{\theta,2}(z) := \kappa_\theta(z^2),
$$
completing the proof for $n = 2$.

Let $n \ge 3$. First, applying the corresponding square functional
$$
J_2(f) = \Bigl(\fc{|c_n(f)|}{M_n}\Bigr)^{2/n},
$$
we derive by the same arguments, as above for the first coefficient
$c_1^0$, that also the second coefficient $c_2^0$ of any extremal
function $f^0$ for $c_n$ must vanish. Hence,
$$
f^0(z) = c_0^0 + c_3^0 z^3 + c_4 z^4 + \dots.
$$

To complete the proof for $n = 3$, one must deal with
the cubic functional
$$
J_3(f) = \fc{|c_3(f)|}{M_3}.
$$
Arguments similar to those applied above provide now
the equalities
$$
|c_3^0| = \kappa_0^\prime(0) = 2/e,
$$
which can hold only for the function $f^0(z) = \kappa_0(z^3)$
and its rotations.

For $n > 3$, comparing successively the relation (3.36) with the functional
$$
J(f)^{n-1} = (|c_n(f)|/M_n)^{(n-1)/n},
$$
one establishes in the same way that the expansion of any
extremal function for $c_n$ must be of the form
$$
f^0(z) = z + c_n^0 z^n + c_{n+1}^0 z^{n+1} + \dots.
$$

So it suffices to deal now with the functions $f \in \mathcal B_1^0$
such that $c_1 = \dots = c_{n-1} = 0$. For such functions, comparison of
the relation (3.36) with the power
$$
J(f)^n = |c_n(f)|/M_n,
$$
provides immediately that the extremal value of $|c_n|$ is
$$
|c_n^0| = \kappa_0^\prime(0) = 2/e.
$$
Therefore, $f^0(z) = \kappa_0(z^n)$, up to rotations.
Theorem 1.1 is proved.

\medskip
{\small\em{
\leftline{Department of Mathematics, Bar-Ilan University}
\leftline{52900 Ramat-Gan, Israel}
\leftline{and Department of Mathematics, University of Virginia,}
\leftline{Charlottesville, VA 22904-4137, USA}}}

\end{document}